\documentclass[conference]{IEEEtran}
\IEEEoverridecommandlockouts
\usepackage{cite}
\usepackage{amsmath,amssymb,amsfonts}
\usepackage{algorithmic}
\usepackage{graphicx}
\usepackage{textcomp}
\usepackage{dsfont}
\usepackage{enumerate}
\usepackage{xcolor}
\def\BibTeX{{\rm B\kern-.05em{\sc i\kern-.025em b}\kern-.08em
    T\kern-.1667em\lower.7ex\hbox{E}\kern-.125emX}}
    
\usepackage{algorithm2e}
 \RestyleAlgo{ruled}    
\usepackage{bbm}
\usepackage{bm}

\newcommand{\TR}[2]{#2}  

\newcommand{\hidespecial}[1]{}  

\newcommand{\bU}{{\bar{\bf U}}} 
\newcommand{\bu}{{\bar{\bf u}}}
\newcommand{\kav}[1]{} 
\newcommand{\sbu}{\bar{U}}
\newcommand{\rsbu}{\bar{u}}

\newcommand{\raU}{\tilde{\bar{u}}}

\renewcommand{\H}{{\bf H}}

\newcommand{\bbeta}{{\bm \beta}}
\newcommand{\gidc}{\hat{\mathds{I}}}

\newcommand{\pof}{\mathds{P}}
\newcommand{\mof}{\mathds{M}}

\newcommand{\lh}{\hat{l}}
\renewcommand{\u}{\hat{U}}
\newcommand{\U}{{\bf U}}  
 
\newcommand{\Y}{{\bf Y}} 
\newcommand{\sU}{U}

\newcommand{\indca}[1]{\mathds{1}\left\{#1\right\}}

\newcommand{\indc}[1]{\mathds{1}_{\left\{#1\right\}}}

\newtheorem{lemma}{Lemma}
\newtheorem{thm}{Theorem}

\newcommand{\eop}{{\hfill $\blacksquare$} }

\newcommand{\N}{{\bf N}}
\newcommand{\C}{{\bf C}}
\newcommand{\X}{{\bf X}}
\newcommand{\Z}{{\bf Z}}

\newcommand{\rx}{{\bf x}}
\newcommand{\bfu}{{\bf u}}
\newcommand{\rc}{{\bf c}}
\newcommand{\rz}{{\bf z}}
\newcommand{\hide}[1]{}
\begin{document}

\title{Fair opportunistic schedulers for Lossy Polling systems  \\
\thanks{The first author's work is partially supported by PMRF, India.}
}

\author{\IEEEauthorblockN{Vartika Singh}
\IEEEauthorblockA{\textit{ IEOR, IIT Bombay, India } \\
vsvartika@iitb.ac.in}
\and
\IEEEauthorblockN{Veeraruna Kavitha}\IEEEauthorblockA{\textit{ IEOR, IIT Bombay, India } \\
vkavitha@iitb.ac.in}}

\maketitle

\begin{abstract}
Polling systems with losses are useful mathematical objects that can model many practical systems like travelling salesman problem with recurrent requests. One of the less studied yet an important aspect in such systems is  the disparity in the utilities derived by the individual stations. Further, the random fluctuations of the travel conditions can have significant impact on the performance. This calls for a scheduler that caters to the fairness aspect, depends upon the travel conditions and the dynamic system state. 

Inspired by the generalized alpha-fair schedulers of wireless networks, we propose a family  of schedulers that further considers binary knowledge of the travel conditions.    These  schedulers are opportunistic, allocate the server    to a station with bad travel condition only when the station has accumulated too little a utility by the decision epoch.  
We illustrate that the disparities among the individual utilities diminish to zero, as fairness factor increases, and further that the price of fairness decreases as the number of stations increase.  

\end{abstract}


\hide{
 \section{paper idea}
 
 \begin{enumerate}
    {\color{purple} \item First part is the conjecture on convergence of $\sbu$.
     \item theorem 1:
     
\begin{enumerate}
    \item for all $\alpha>0$, there is some $\epsilon$, such that $\max |\rsbu_i- \rsbu_j|<\epsilon_\alpha$
    \item $\epsilon_\alpha$ decreases with $\alpha$
    \item  $\epsilon_\alpha$ decreases to zero
    \item {\color{blue} Done!}
\end{enumerate}
\item theorem 2 (PoF) 
: If $\{\lambda\}_i$ are all equal, $s_i$ are all equal, if travel condition on good and bad is identical, then
\begin{enumerate}
    \item $0<g_{ij}=P(T_{ij}=G)<1$ $\forall$ $i,j$
    \item Pof$(\alpha=\infty)\to 0$ as $n\to \infty$
    \item Pof$(\alpha)\to$ Pof$_g(\alpha)$ (strongest, implies b) for symmetric conditions / may try for asymmetric conditions
\end{enumerate}}
\item the lemma that server does not stay in any queue infinitely
\item utilities start closing in with more $\alpha$, even for different $\lambda$
\item the utilities at $\alpha=0$ itself start closing in, as $m\to \infty$ for same $\lambda$
\item if we have time, we can try to show the pof converges to the pof of good system (will require the maximum theorem kind of things)
\item we are assuming $\rsbu_i\ge 0$, can we prove / justify it
\item opportunistic scheduler has to be explained in introduction
 \end{enumerate}
 }

\vspace{-5mm}
\section{Introduction}

The motivation for this study comes from a variant of travelling salesman problem, where the salesman keeps travelling from one location to another to provide some service. It is natural that some locations will be  preferred more by the salesman in comparison to others, possibly because of   travel conditions, demands, etc. This may lead to the starvation of less preferred locations and calls for a fair policy.

We consider a single server that visits several stations  to serve the demands. After each service, it makes a decision to either remain in the same station, or  to travel to other stations  depending upon travel conditions, system state, etc.
In all, we have a polling system (e.g., \cite{polling},\cite{PollSchedule}, \cite{fairpoll})  with controlled (dynamic) service and random switching times. 



Fairness is a relatively  less studied aspect in polling systems (e.g., \cite{fairpoll}) and our focus is on utilities of individual stations and further on opportunistic aspects.
We take inspiration from generalised $\alpha$-fair opportunistic scheduler of the wireless networks (e.g., \cite{cellular,khushnerfair})  which considers fair allocation of  resources (like bandwidth) among users over  time; authors in \cite{mayurNH} argue that the price of fairness (PoF) of opportunistic schedulers for wireless networks with large users  is negligible. We propose a Fair opportunistic Polling Scheduler (FoPS) that prescribes the next station to be visited by the server, depending upon current server location, travel conditions and number of customers waiting at different stations; to cater to fairness aspects, FoPS also depends upon the accumulated utilities of the individual stations till the decision epoch. The FoPS  maintains a balance between the efficiency and the fairness dictated by a single parameter $\alpha$. Basically, it chooses  stations with good travel conditions as long as none of the stations are starved.

The salient features and contributions of our work are:  i) we propose a consolidated objective function that considers losses (occurred by losing the customers) and the rewards (obtained by serving customers);   ii) the proposed scheduler is based on some anticipated utilities of the individual stations that depend upon the current state and travel conditions; iii) we define a measure of fairness (MoF) that quantifies the disparity in utilities of all the stations and show that the metric reduces with increase in $\alpha$ or the number of stations; and   finally,  iv) the PoF, i.e., the price of fairness (a metric similar to price of anarchy) defined in  \cite{Bertsimas,mayurNH}, is miraculously negligible even with heterogeneous demands. 

To summarize, the proposed schedulers achieve various levels of fairness (the most fair policy results in  almost equal individual station-utilities), with almost negligible loss of efficiency in several scenarios.   
In other scenarios, the price of fairness reduces to a value that can be achieved in an ideal  system, where the travel conditions are always good and this is possible with large number of stations. 
Surprisingly,  the efficient scheduler among FoPS is itself fair with large number of stations. We perform numerical simulations to reaffirm our findings, and observe around ten stations suffice for negligible MoF, this number further decreases with increased $\alpha$. Interestingly, PoF is negligible in all the case studies.


\hide{Our contributions are as follows:  i) we propose a consolidated objective function that considers losses incurred by loosing the customers due to finite buffer length,  and the rewards obtained by serving the  customers; ii) we propose a fair opportunistic polling scheduler, inspired by wireless opportunistic schedulers, to  maintain fairness in the system;   iii) further,  the proposed scheduler also considers travel conditions  along with other state components; iv) we show that the price of fairness, a metric similar to price of anarchy defined in (\cite{Berskas,M}), is miraculously negligible.

We observe that the fairness is achieved (all stations achieve almost equal individual utilities), with almost negligible loss of efficiency in several conditions.   
In the other conditions, the price of fairness reduces to a value that can be achieved without disparity in channel conditions, where all the channels are in their best condition, and this is possible with large number of stations. Surprisingly,  the efficient scheduler among fair opportunistic polling schedulers attains fairness with large number of stations.}

Polling systems are well studied objects in literature (see \cite{polling} for an elaborate and recent survey) and we discuss a few relevant strands here; 
  authors in \cite{fairpoll} study fairness aspects with respect to various service disciplines (within queues) like exhaustive, gated,  etc.;  \cite{PollSchedule} considers dynamic scheduling policies 
   to schedule jobs within each queue.
However, to the best of our knowledge, none of these papers   consider fairness aspects across stations and opportunistic advantage resulting with some minimal knowledge about switching conditions. Further, they do not discuss the price and measure of fairness.

\vspace{-1mm}
\section {Problem statement and background}
There are $m$ stations and a single server. At any station $i$, the demands  arrive according to an independent Poisson process with rate $\lambda_i$  and join the queue,  $Q_i$. The server travels from one station to another and keeps serving the waiting customers. The service times are independent and identically distributed (i.i.d.) across all the stations with mean $1/\mu_s$. Further the buffer length at station $i$ is  $b_i$, i.e., if a new demand arrives at station $i$, and if the queue length is $b_i$, such a demand is lost.

We consider a controlled service, in which the server upon arrival to station $i$ serves exactly one customer (if the queue is non-empty); it then makes a  decision either to stay in station $i$ or to travel to some other station, depending upon on the system state. If the queue $Q_i$ is empty, the server waits for some time and then makes a new decision; for mathematical simplicity, we model the waiting time distribution to be the  same as the service time distribution. The server does not serve any demand during wait times. The state of the system is primarily described by the number of customers waiting at different stations and the travel conditions.
From any station to any other station, the travel conditions can either be good  or bad. When the travel conditions are good (respectively bad), the travel time is distributed according to a known distribution (same for all stations) with mean $\mu_g$ (respectively $\mu_b$); naturally $\mu_g<  \mu_b$. {The travel condition from station $i$ to $j$ is good with probability $0<p_{ij}<1$.}

\hide{The travelling time from any station to any other station has some known distribution with finite mean depending upon the travel condition. There are two possibilities,
\begin{enumerate}
    \item Good travel condition with probability $p_i$,
    \item Bad travel condition with probability $(1-p_i)$.
\end{enumerate}If the travel conditions are good, the travel time is distributed (same for all the station pairs) with a smaller mean and variance $(\mu_g,\sigma_g)$, and if the travel conditions are bad, the travel time has a larger mean and variance $(\mu_b,\sigma_b)$ (same for all the station pairs). Given the same travel conditions, the travel time is distributed according to same known distribution for all station pairs.}

In all, we have a polling system (e.g., \cite{polling}) with controlled (dynamic) service and random switching times. Further, the controlled dynamic decision also depends upon the travel conditions to the rest of the stations; relevant travel conditions are assumed to be known at the decision epochs. In our model, depending upon the controlled decision, the server can be idle and wait for next arrival in the queue even if there are demands at other stations. Furthermore, we consider fairness aspect, that attempts to reduce the disparities in individual station  utilities.



Fairness is a well studied concept in the context of wireless networks (see \cite{mayurNH}, \cite{priority} etc.). Our approach is to adopt the idea from wireless networks to design a fair policy. 
Before we discuss these details, we digress briefly to discuss the fairness concepts in wireless networks.

\vspace{-4mm}
\subsection{Motivation: Wireless $\alpha$-fair scheduler}
In a wireless network, a scheduler allocates the resource to one of the users in each time-slot/epoch, and aims to maximize the sum of utilities of all the users. The utilities of any given user are i.i.d  across epochs, but, different users can have different distributions. The scheduler can observe the utilities before allocating the resource.
To maximize\footnote{
\label{foot_one}
Recall the utilities are independent across the epochs and are identical for the same user. Thus efficient scheduler maximizes sum of the expected utilities of all the users, which by Ergodic theorem, almost surely equals the sum of the time averages. As we will see soon, in our context the utilities are not independent across epochs; thus to optimize the sum of the expected utilities in strict sense, one requires Markov Decision Process (MDP) framework. This would be an interesting future aspect that we would investigate.} the sum of utilities, an efficient scheduler allocates the resource to the one with the highest utility in each epoch. This leads to starvation of some of the users  over time, especially ones with lower expected utilities. In such scenarios, it becomes important to be fair while deriving the optimal resource allocation policy.

Generalized $\alpha$-fair schedulers maximize a concave function of time-accumulated utilities (by Ergodic theorem, expected utilities) of all the users, and achieve various levels of fairness depending upon parameter $\alpha$ (e.g., \cite{cellular,khushnerfair}):

\vspace{-4mm}
{\small\begin{eqnarray}\label{eqn_fair_concave}
 \sup_{\bbeta = (\beta_1, \cdots, \beta_m)} \sum_i  \Gamma_\alpha \left ( \rsbu_i(\bbeta) \right ) \text{ with } \rsbu_i(\bbeta):=E[ \sU_{i} \beta_i (\U)]\ ],\hspace{-10mm} \nonumber  \\
 \Gamma_\alpha(\rsbu) := \frac{\rsbu^{1-\alpha} \indc{   \alpha \neq 1 } }{1-\alpha} + \log(\rsbu) \indc{  \alpha =1  },
 \label{Eqn_alpha_fair}
\end{eqnarray}}where $\U := (U_1, \cdots U_m)$.
In the above, max-min fairness is achieved with $\alpha\to\infty$, proportional fairness is achieved with $\alpha=1$, and $\alpha=0$ achieves efficient scheduler.

The solution $\bbeta^{(\alpha)}   = (\beta_1, \dots, \beta_m)$ of \eqref{eqn_fair_concave}, for any $\alpha$, is given by the solution $\bu^{(\alpha)}=(\rsbu_1,\dots, \rsbu_m)$ of following $m$-dimensional fixed point equation  (\cite{cellular,khushnerfair}):

\vspace{-4mm}
{\small \begin{eqnarray}
\label{eqn_beta_alpha_FP}
\rsbu_n = E[ \sU_n \beta_n (\U) ], \     \beta_n(\U) :=\prod_{i \ne n}  \indca{   \frac{\sU_n}{ \left (\rsbu_n \right )^\alpha }     \ge  \frac{\sU_i}{ \left (\rsbu_i \right )^\alpha } }.\hspace{-1mm} \hspace{-1mm}
\end{eqnarray}}
An
   iterative algorithm  that uses   the average utilities  derived by the individual users till epoch $k$, represented by $\bU_k = (\sbu_{1,k},\cdots,\sbu_{m,k})$,   is proposed to attain the above solution:
   
   \vspace{-4mm}
   {\small
\begin{eqnarray}\label{eqn_iterative}
\sbu_{n,k+1} &=& \sbu_{n,k} + \frac{1}{k+1} \left ( U_{n,k+1} \beta_{n,k+1} - \sbu_{n,k} \right ),  \nonumber \\
\beta_{n,k+1} &=& \prod_{i \ne n}  \indca {  \frac{U_{n,k+1}}{ \left (\sbu_{n,k} \right )^\alpha }     \ge  \frac{U_{i,k+1}}{ \left (\sbu_{i,k} \right )^\alpha }  }.  \label{Eqn_beta_scheduler}
\end{eqnarray}}
In \cite{khushnerfair,cellular}, it is proved that the above algorithm converges weakly to asymptotic utilities \eqref{eqn_beta_alpha_FP} for i.i.d. channels. 
 Thus the data scheduler  $\bbeta_{k+1}$, for time slot $k$,  depends on the channel estimates $\U_{k+1}$ as well as the average utilities 
 $\bU_{k}$.

Further, there must be a price of fairness for deviating from the efficient scheduler. But the results  in \cite{mayurNH} show that the price of fairness is negligible if there are a large number of users. We would also investigate this aspect. 

\section{Fair Dynamic Opportunistic schedulers}
The idea is to derive a fair  policy for  polling systems following the guidelines from wireless networks. In our case, server is the resource to be allocated to the users (stations).  We need to define the instantaneous utilities of the stations in such a way that, the utility of the server is the sum of instantaneous utilities of all the stations over time. The instantaneous utility at any station $i$ should consist of gains from serving the customers at $i$ and losses due to customer drops because of buffer length $b_i$. We begin with some definitions.

Let $\N_k = (N_{1,k}, \cdots, N_{m,k})$  be the number of customers waiting at various stations, and $S_{k-1} \in \{1,\cdots,m\}$ be the server location at the beginning of $k$-th decision epoch. 
The flag  $C_{i,k}\in\{g,b\}$ represents the travel condition, if the server has to travel from  $S_{k-1}$ to $i$. Set $C_{i,k} = 0$ if $i = S_{k-1}$, to indicate no travel is required, and let 
$\C_k = (C_{1,k},\cdots, C_{m,k}) $.


The scheduler in \eqref{eqn_iterative} depends upon the average utilities $ \bU_{k} $ derived by individuals till the decision epoch $k$. Observe from \eqref{eqn_iterative} that the users with smaller $\sbu_{i,k}$ have higher weight factor; higher  $\alpha$ magnifies to a larger extent. \textit{The scheduler also considers opportunistic aspect as the scheduling decision additionally depends upon the instantaneous utilities $\U_k = \{U_{i,k}\}_i$}. Our aim is similar, average utilities as defined below can again be used for similar purpose, opportunities provided by the travel conditions $\C_k$ and the current state, $(S_{k-1}, \N_k)$ can also be used:  first define the average utility obtained till the decision epoch $k$ for any $Q_i$ (factor $\gamma \approx 1$ helps in mathematical tractability), 

\vspace{-3mm}
\small{\begin{equation}\label{eqn_u_bar_average}
%
\sbu_{i,k} = \sbu_{i,k-1} + \frac{1}{k^\gamma} \left ( {\u}_{i,k} - \sbu_{i,k-1} \right ),
\end{equation}}where ${\u}_{i,k}$ has to be appropriately defined to reflect the opportunities provided by  $\Z_k:=(\X_k,\C_k)$ with $\X_k= (S_{k-1},\N_k, \bU_k)$; we will soon see that $\u_{i,k}$ will be anticipated utilities of some  $U_{i,k}$, the actual instantaneous utilities gained by the decision taken at epoch $k$; we begin with defining $U_{i,k}$ while $\u_{i,k}$ are defined in equation \eqref{eqn_util}. 

As already mentioned, the  utility of any  station is comprised of gains and losses and depends upon the scheduling decision. Say ${\bbeta}_k= \bbeta(\Z_k)=(\beta_{1,k},\dots,\beta_{m,k})$ is the decision, where $\beta_{j,k}$ is the probability of travelling from $Q_{S_{k-1}}$ to $Q_j$, with  $\beta_{j,k}\in[0,1]$ and $\sum_{j} \beta_{j,k}=1$. \textit{Station $i$ gains a reward $w$ if one of it's users is served}. It's losses are captured by the number of users that have overflown during the travel time (if any) and the service time. Thus the instantaneous utility of station $i$ because of  decision $\bbeta_k$ equals,

\vspace{-5mm}
{\small
\begin{eqnarray}\label{eqn_actual_util}
{U}_{i,k}(\bbeta_k,\Z_k)=w\indc{S_k=i}\indc{{\tilde{N}}_{i,k}\ge 1} - \sum_{j}\indc{S_{k}=j}{L}_i(S_{k-1},j), \hspace{-2mm}
\end{eqnarray}}%
where $\tilde{N}_{i,k}$ is number of customers waiting in $Q_i$ just before the service starts (can be different from $N_{i,k}$ if the server travels), $L_i$ represents the losses in the same queue $Q_i$ during the time interval between the two decision epochs and probability  $P(S_{k}=j|\bbeta_k)=\beta_{j,k}$ for any $j$. The losses $L_i$ depend upon the travel conditions, service times, and hence upon $S_{k-1}$ and the next station chosen.

We immediately have the following observation. The state $\{\Z_k\}$ (in fact $\{\X_k\}$) evolves as a non-homogeneous Markov Chain with countable state space under any given dynamic policy $\bbeta$ that depends only upon the system state $\Z$.
\vspace{-1mm}
\subsection{Opportunistic utilities based on anticipation}
It is clear that the utilities in \eqref{eqn_actual_util} are not known at the decision epoch $k$. Hence we propose to use the expected values of $\{U_{i,k}\}_i$, conditioned on the available information $\Z_k$. These form instantaneous utilities for $k$-th decision epoch and equal,
%

\vspace{-4mm}
{\small
\begin{eqnarray}\label{eqn_util}
    \u_{i,k} (\bbeta,\Z_k) = \beta_{i,k}  \gidc_i(\Z_k)  - \sum_{j  } \beta_{j,k}  \lh_i(j,\Z_k),
\end{eqnarray}}where i)
$\gidc_i (\Z_k):= w E[ \tilde{N}_{i,k} > 0 | \Z_k, S_k = i] $ equals  the expected gain at $Q_i$ by serving a customer; and ii) ${\lh}_i (j,\Z_k):=E[ L_i(S_{k-1},S_k)| \Z_k, S_k = j] $ is the expected loss at $Q_i$ (due to buffer length $b_i$)  when station $j$ is chosen for next service.

By conditioning on $\Z_k=\rz = (\rx,\rc)$, the expected losses and gains of station $i$ are given by,

\vspace{-4mm}
{\small \begin{eqnarray}\label{eqn_gains}
\gidc_i (\rz)  &=&  w \left (\indc{n_{i}>0} + \indc{n_{i}=0, s\ne i} E\left[\indc{{\cal N}_i (T_{s,i}   )  > 0 }| c_{i}\right] \right ),  \nonumber \\ 
{ \lh}_i (j,\rz)  &=&  E \left [\left ({\cal N}_i (T_{s,j} + J ) + n_{i} - b_i \right )^+  \left | c_{j}, n_{i}  \right .  \right  ],\hspace{1mm}
\end{eqnarray}}where ${\cal N}_i (T_{s,i})  $ is the number of Poisson arrivals in $Q_i$ during the travel time $T_{s,i}$ and $J$ is the service time (waiting time if queue is empty). %
Further the anticipated losses of any station  are  bounded uniformly by $\lh^* := \max_q \{\lambda_q (\mu_b + 1/\mu_s)\}$:

\vspace{-4mm}
{\small 
\begin{eqnarray}\label{eqn_loss_bound}
\lh_i(j,\rz) &=& E [\left ({\cal N}_i (T_{s,j} +  J ) + n_{i} - b_i \right )^+  | c_{j}, n_{i}    ], \\
&&\hspace{-5mm}\le \  E [{\cal N}_i (T_{s,j} + J)|c_j]
\TR{=E[\lambda_i (T_{s,j} + J)| T_{s,j}, J],\nonumber\\ &&\hspace{-5mm}\le}{\le} 
\lambda_i (\mu_b + 1/\mu_s) \le \lh^* .
\nonumber
\end{eqnarray}}

\subsection{Fair opportunistic Polling Scheduler (FoPS($\alpha$))}

We are now ready to describe the proposed scheduler.
Taking inspiration from wireless scheduler given in \eqref{eqn_beta_alpha_FP} and \eqref{eqn_iterative}, we propose the following scheduler, parameterized by $\alpha$, for decision epoch $k$ when state $\Z_k=\rz=(\rx,\rc)$,
\begin{equation}\label{eqn_dyn_dec}
\mbox{\bf FoPS}(\alpha):  \hspace{10mm}  \bbeta_k^*(\rz) =  \arg\max_\bbeta \sum_i   \frac{\u_{i,k}(\bbeta, \rz)}{\raU_{i}^\alpha}, \hspace{4mm}
 \end{equation}where $\raU:=\max\{\delta,\rsbu_{i}\}$; here  $\delta>0$ is a small quantity  introduced to take care of possible negative values. Observe that this scheduler is lot more complicated than the scheduler in \eqref{eqn_iterative}: a) in wireless networks, if a user is not chosen at an epoch, it's instantaneous utility is zero, however in our system the unselected stations get negative utility; b) further the negative utility depends upon travel condition of the chosen path (station); c) thus the utility of any station in any decision epoch just does not depend upon whether the station is selected or not; d) and hence the system evolves according to a non-homogeneous Markov chain driven by dynamic decision \eqref{eqn_dyn_dec}.


We can re-write the decision \eqref{eqn_dyn_dec}  as,

\vspace{-4mm}
{\small
\begin{eqnarray}\label{eqn_ob}
\bbeta_k(\rz)& =&  \arg\max_{\tilde{\bbeta}} \sum_i \tilde{\beta_i} O_i(\rz)\ \  \mbox{where,} \nonumber\\
O_i(\rz) &:=& \frac{{\gidc}_i(\rz)}{\raU_{i}^\alpha} - \sum_j\frac{ { \lh}_j(i,\rz)}{\raU_{j}^\alpha}.
\end{eqnarray}}In the above expression, $O_i$ can be interpreted as the fair-weighted utility of the server upon choosing $Q_i$. Hence,   the decision $\bbeta_k(\rz)$ is given by,
\begin{equation}\label{eqn_choice_station}
\beta_{q,k} (\rz)=\frac{\indc{q \in \arg\max_i O_i(\rz)}}{\left|\arg\max_i O_i(\rz)\right|}, \mbox{ for any } q.
\end{equation}

\hide{\color{red} Given any SMR policy $\beta$, by Ergodic theorem ??
\begin{align}
&\lim_{T\to \infty} \frac{1}{T} \sum_{k=1}^T \sum_{n=1}^m  \beta_{n} (\Z_k) \sum_j {U}_j(n)=\nonumber\\ &E^{\pi(\beta)}_{x_0}\left[\sum_{n=1}^m \beta_n(\Z)\sum_j {U}_j (n) \right] \ a.s. := \sum_j \rsbu_j(\beta,x_0)
\end{align}

Writing it user-wise for any user $j$:
\begin{eqnarray}
\lim_{T\to \infty} \frac{1}{T} \sum_{k=1}^T \sum_{n=1}^m  \beta_{n} (\X_k) { U}_j(K,n,\X_k) &=&\\ E^{\pi(\beta)}_{x_0}\left[\sum_{n=1}^m \beta_n(\X) {U}_j (K,n,\X) \right] \ a.s. :=  \bu_j(\beta,x_0)
\end{eqnarray}
$${U}_j(K,n,\X)= \indc{n=j}\indc{\Tilde{N}_j \ge 1} - \Tilde{L}_j,$$where $\Tilde{N}_j$ is the actual number in queue $j$, just before the service starts, whose conditional expectation was $\bar{p}_j(K,j,\X)$ and $\Tilde{L}_j$ is the actual losses in the queue $j$ just after service/waiting time completion at $Q_n$ whose conditional expectation was $\bar{l}_j(K,n,\X)$
}

We conclude this section with a conjecture that describes the limits of the time-averages  in \eqref{eqn_u_bar_average}, which is required for Theorem \ref{thm_bdd}. Further the conjecture provides interesting and important insights which could be of independent interest. The conjecture is verified using several numerical examples (e.g., Figure \ref{fig:trajectories}) and is partially supported by Theorem \ref{Thm_conv} which also proves the convergence of utilities of \eqref{eqn_u_bar_average}. Towards this we require some important definitions.

\hide{
\begin{figure}
\vspace{-10mm}
    \centering
    \includegraphics[scale=.25]{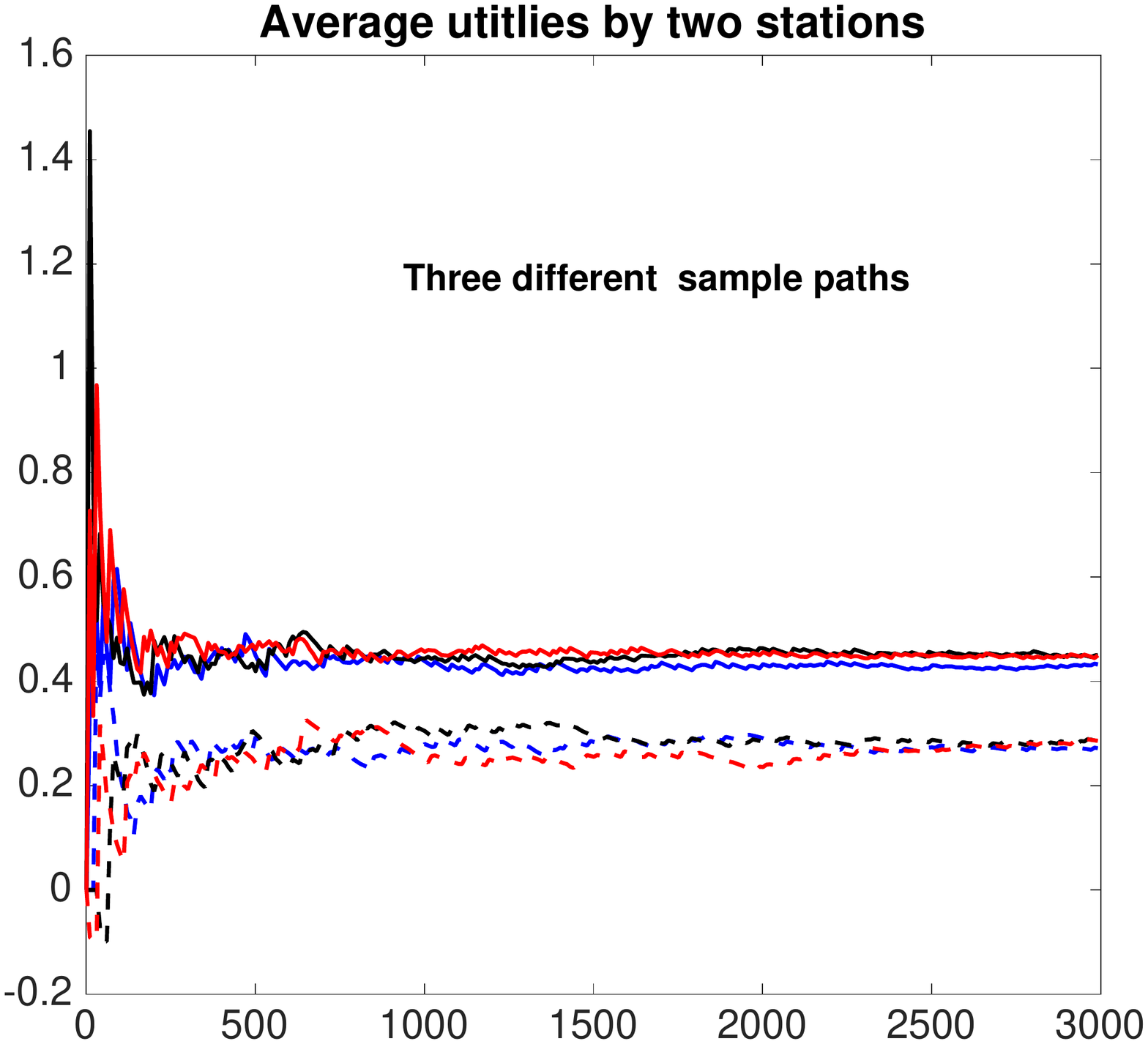}
    \vspace{-15mm}
    \caption{All three sample paths converge to the same limit}
    \label{fig:trajectories}
\end{figure}}

Let $\Y_k(\bu):=(S_{k-1},\N_k,\C_k)$ be the  chain whose transitions are similar to the transitions of corresponding components of $\Z_k$ except that $\bbeta_k$ in \eqref{eqn_choice_station} is now defined using fixed vector $\bu$. {\it Observe that $\mathds{Y}(\bu) =  \{\Y_k(\bu)\}_k$ is a finite-state Markov chain, can be reducible or irreducible depending upon $\bu$} and will have unique stationary distribution in each irreducible closed class for any $\bu$. Further, it satisfies the standard Ergodic theorem (e.g., \cite{LLN}) for any $\bu$.  
These Ergodic chains, in our opinion, form the basis to derive the almost sure (a.s.) limit of time-average utilities defined in \eqref{eqn_u_bar_average} under an assumption that uses the following definition.

{
\begin{enumerate}[{\bf D}.0]
    \item The vector $\bu^*$ is said to be Ergodic-interior if there exists an $\epsilon>0$, such that $\mathds{Y}(\bu)$ has transition probabilities same as $\mathds{Y}(\bu^*)$ for all  $\bu \in \{\bu': |\bu' - \bu^*| < \epsilon\}$.
\end{enumerate}
}

\newcommand{\ta}{ {\bar V}}
\newcommand{\vta}{{\bf   {\bar v}}}
\newcommand{\rta}{{  {\bar v}}}

\begin{thm}{\bf [Ergodicity]}
\label{Thm_conv}
Assume $\gamma > 1$.
For any $\alpha$ and $i$, the time-average utilities,  \vspace{-6mm}
{\small
\begin{equation}\label{eqn_u_bar_act_average} \hspace{20mm}
\ta_{i,k} := \frac{1}{k} \sum_{t=1}^k \sU_{i, t},
\end{equation}}
and the utilities defined in \eqref{eqn_u_bar_average} 
converge. 
The respective limits $\vta, \bu $  can depend upon sample path $\omega$, and further:
\begin{eqnarray}\label{eqn_fixed_point_alt}
\rta_{i} &=& E\left[ \u_i(\bbeta,\Z) \right]
\mbox{ for any $i$,}
\end{eqnarray}
where the expectation is under a stationary distribution (SD) of the Ergodic Markov chain $\mathds{Y}(\bu)$, if 
$\bu$ is Ergodic-interior as defined in {\bf D}.0.
The transitions of Markov chain $\mathds{Y}(\bu)$ are captured by:

\vspace{-4mm}
{\small \begin{eqnarray}\label{eqn_fixed_point_2}
 \beta_q(\Y) &=&\frac{\indc{q \in \arg\max_i O_i(\Y, \bu)}}{\left|\arg\max_i O_i(\Y, \bu)\right|}, 
\mbox{ with } \\
O_i(\Y,\bu)  &:=&  \frac{{\gidc}_i(\Y, \bu)}{\raU_{i}^\alpha} -\hspace{-1mm} \sum_j\frac{ { \lh}_j(i,\Y, \bu)}{\raU_{j}^\alpha}, \ \raU_i:=\max\{\delta,\rsbu_{i}\},\nonumber\\
   \u_{i} (\bbeta,\Y, \bu) &=& \beta_{i}  \gidc_i(\Y, \bu)  - \sum_{j  } \beta_{j}  \lh_i(j,\Y, \bu).\nonumber 
 \end{eqnarray}} 
 
\end{thm}
\textbf{Proof} is in Appendix \eop

{\bf A conjecture:}
When $\gamma $ is close to $1$, we anticipate the two limits of the above theorem $\bu$ and $\vta$ to be close (observe $\vta$ also equals the time limit   of \eqref{eqn_u_bar_average} but with $\gamma = 1$, as ${\bf \u}$ are conditional expectations of $\U$). For now we conjecture them to be approximately equal,  further using \eqref{eqn_fixed_point_alt} conjecture that the following fixed point equation is satisfied as $\gamma \downarrow 1$,\vspace{-4mm}
\begin{eqnarray}\label{eqn_fixed_point}
\rsbu_{i} &=& E\left[ \u_i(\bbeta,\Z) \right]
\mbox{ for any $i$.}
\end{eqnarray}

{\bf Remarks:} i) The  proof of the conjecture is immediate for $\alpha=0$ by the well-known Ergodic theorem provided in   standard text books (see also \cite{LLN}) when $\gamma = 1$; $\{\Z_k\}$ is then finite state space irreducible Markov chain from \eqref{eqn_ob} under {\bf A.1} (given in  section \ref{sec_fair_analysis}). ii) With $\alpha>0$, $\{\Z_k\}$ is a non-homogeneous Markov chain and we are working towards the proof  of the conjecture for this case. For this short paper, we continue with further analysis by assuming the conjecture to be true. 
%
 iii) If the conjecture were true, 
the average utilities of \eqref{eqn_u_bar_average} converge to the limit utilities $\bu = \bu^{(\alpha)}$, which are given by  stationary expected values of an appropriate Markov chain $\mathds{Y}$ with `limit transition probabilities' given in \eqref{eqn_fixed_point_2}; these are governed by fixed point equations \eqref{eqn_fixed_point} as the limit transition probabilities depend upon the limit utilities $\bu$.
{\it iv) the   conjecture  is used  only  in the proof of Theorem \ref{thm_bdd}, while the rest of the results are derived independently.}

\begin{figure*}   
\begin{minipage}{0.14\textwidth}
\vspace{-5mm}
    \includegraphics[width=3.2cm,height=3.4cm]{Trajectories.pdf}
    \vspace{-13mm}
    \caption{All three sample paths converge to the same limit}
    \label{fig:trajectories}
    \end{minipage}
    \hspace{10mm}
\begin{minipage}{0.35\textwidth}
\begin{minipage}{0.4 \textwidth}
    \centering
     \includegraphics[scale=0.25]{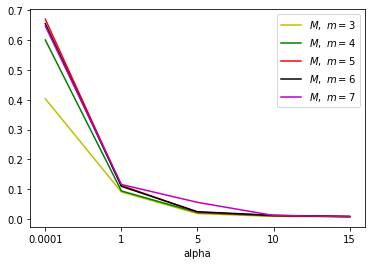} 
\end{minipage}\hspace{10mm}%
\begin{minipage}{0.4 \textwidth}
    \centering
     \includegraphics[scale=0.25]{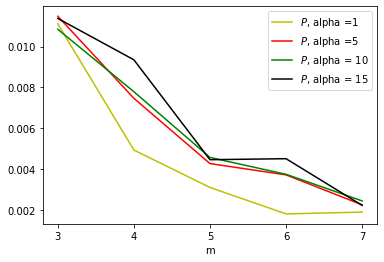} 
    
\end{minipage}
\caption{{\bf MoF $\mof$ and PoF $\pof$ versus $\alpha $ and $m$}: $\lambda=0.67/m$, $\mu_g=2$, $\mu_b=6$, $\sigma_g = \sigma_b = 0.1$, $b=5$, $\mu_s =1/3$, $p_{ij}\in\{0.1,0.9\}$, $w=6$.
\label{fig:mof_vs_alpha}}
\end{minipage}
\hspace{10mm}
    \begin{minipage}{0.36\textwidth}
\vspace{-4mm}
    \centering
    \hspace{-16mm}
    \begin{minipage}{0.3\textwidth}
         \includegraphics[width=3.2cm,height=3.4cm]{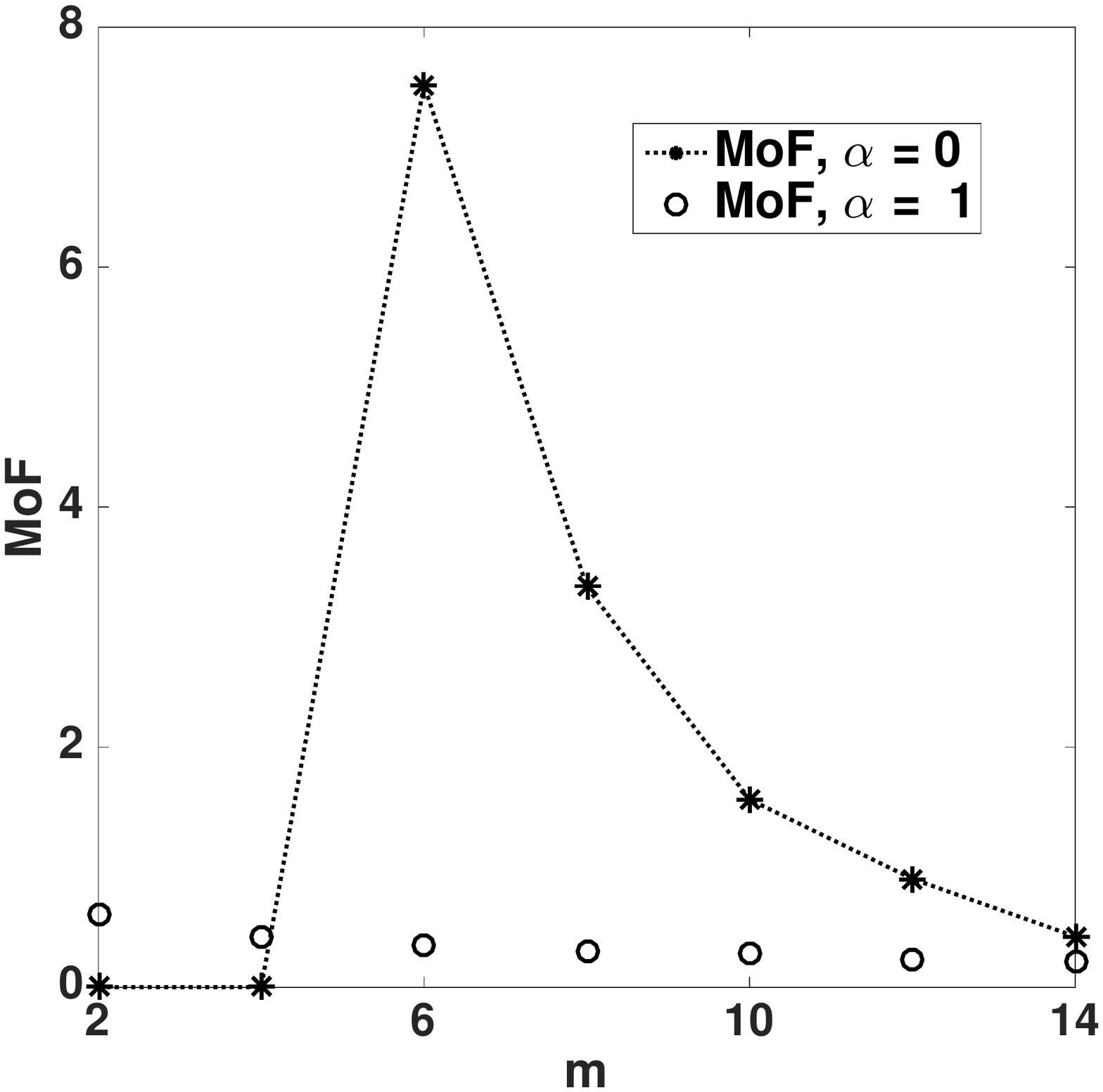}
    \end{minipage}
         \hspace{10mm}
        \begin{minipage}{0.3\textwidth}
        \vspace{-2mm}
         \includegraphics[width=3.2cm,height=3.2cm]{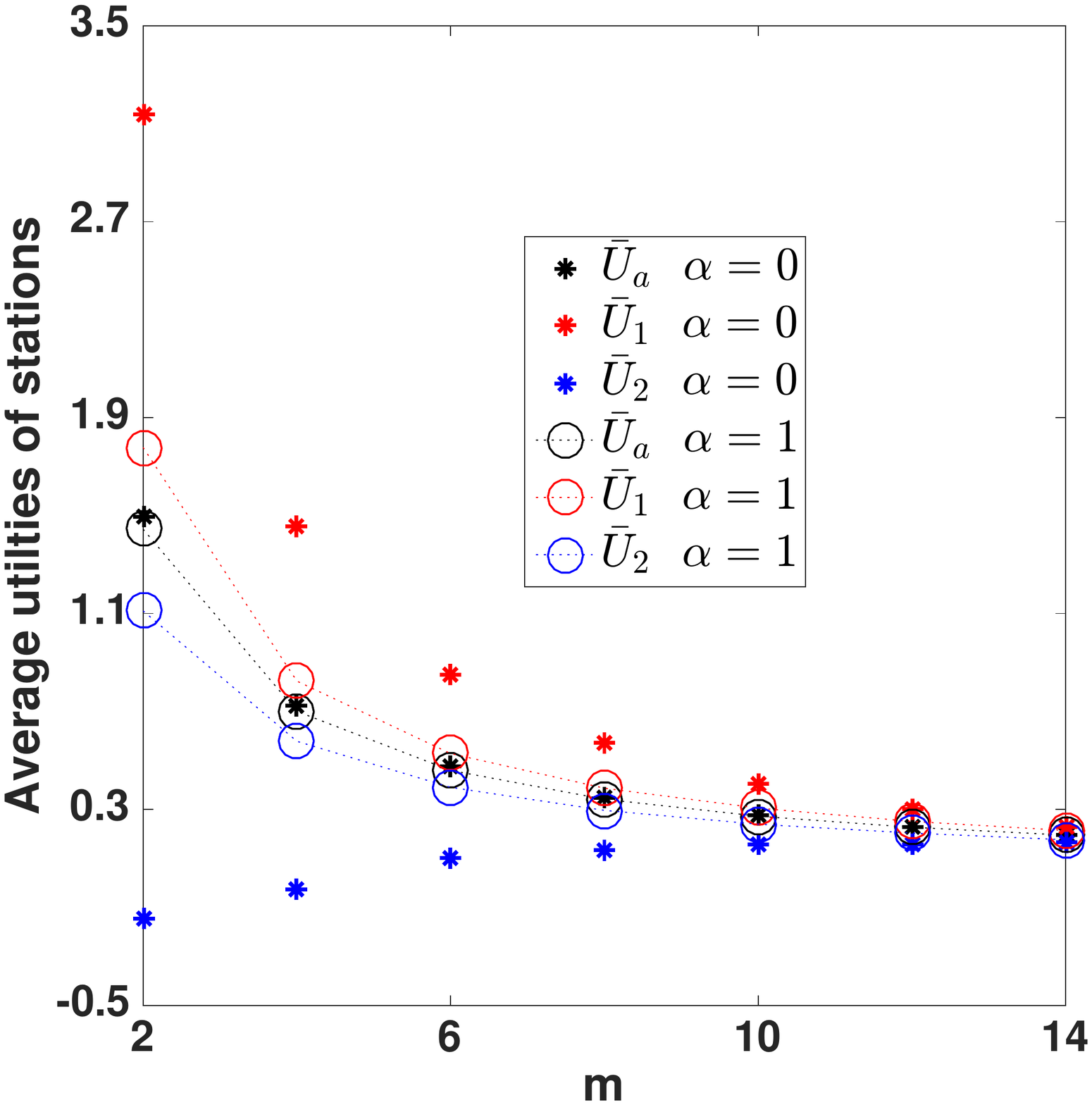}
    \end{minipage}
        \vspace{-8mm}
    \caption{
{\bf Good stations with more demands:}   
    $\mu_g = 2$, 
$\sigma_g = 0.01$,
$\mu_b = 6$
$\sigma_b = .08$,
$b_i \equiv 4$, $\mu_s = 1/3$ $p_{i,j} \in \{ 0.1, 0.9 \}, \lambda_i \in \{0.37 /m, 0.93 /m\}$ and $w = 4$ 
    \label{fig:mof_and_util}}
    \end{minipage}
    \vspace{-4mm}
\end{figure*}

\section{Fairness Analysis} \label{sec_fair_analysis}
Our aim in this section is to derive the analysis of the proposed FoPS($\alpha$) schedulers. We begin with defining two measures that respectively capture the `amount of fairness' achieved and the `price paid for that fairness'. 

A scheduler is said to be efficient if it maximizes the system/server utility $\sum_i \rsbu_i$. From \eqref{eqn_dyn_dec}, FoPS$(0)$ is efficient among the FoPS($\alpha$) schedulers. The server utility is maximized by efficient scheduler, but the individual utilities of the stations $\{\rsbu_i\}$ can have \textit{huge disparities depending upon the travel condition statistics and the arrival rates of individual stations. The precise aim of the fair schedulers is  to reduce this disparity.} There is an anticipated reduction in the server utility and hence efficiency, when one seeks fairness. The price of fairness (PoF) of any fair scheduler captures this price as,
\begin{equation}\label{eqn_def_pof}
  \pof(\alpha):=  \frac{\sum_{j=1}^m \rsbu_j(\alpha)-\sum_{j=1}^m\rsbu_j( 0)}{\sum_{j=1}^m\rsbu_j(0)},
\end{equation}
where $\rsbu_j(\alpha)$ is solution of the fixed point equation \eqref{eqn_fixed_point} (or the limit of \eqref{eqn_u_bar_average}) for given $\alpha$.
One can say a scheduler is max-min fair (see \cite{cellular} and references therein) if it ensures that the utilities of all the stations are equal. Hence we define the measure of fairness (MoF) of any scheduler using  the amount of disparities between the individual utilities as below: \vspace{-4mm}
\begin{equation}\label{eqn_def_mof}
    \mof(\alpha) = \max_{i,j} \frac{\raU_i(\alpha) }{\raU_j(\alpha) }-1.
\end{equation}
By definition  $\mof \ge 0$,  
the smaller the value, the `fairer' is the system. Further, $\mof$ can also be zero if all the utilities $\rsbu_i < \delta$. We make the following assumption for further analysis,
\begin{enumerate}[{\bf A}.1]
    \item Maximum possible losses should be less than maximum possible gain, i.e., $$\rho_B := \max_{i,q}  \sum_j {\cal N}_j E [( T_{i,q}+J)\left |c_q=b \right . ]<w.$$
    
\end{enumerate}

The first   observation about the FoPSs is that  all the stations are served \textit{infinitely often (i.o.)}, i.e., the number of visits to any station increases to $\infty$ (proof is in Appendix): 
\begin{lemma}\label{lem_infinite_stay}
Define the event,
  {\small${\cal I}:= \cap_i \{S_t= i \mbox{ i.o.}\} $}, where $\{S_t= i \mbox{ i.o.}\}=    \cap_t \cup_{k\ge t} \{S_k= i  \} $ and  assume {\bf A}.1. Then for any $\alpha>0$,  we have $P({\cal I})=1$.  \eop
\end{lemma}
%
Thus, with $\alpha>0$, all the stations are served infinitely many times, however the more important aspect is the measure of fairness. Towards this, we have the next result, which shows that the measure of fairness reduces, i.e., `fairness' increases  as $\alpha$ increases. The proof (in Appendix) of this theorem uses the Conjecture \ref{Thm_conv}, and hence is true subjected to the verification of the latter. Nonetheless, the simulations in the next section illustrate that the $\mof$ indeed decreases with an increase in $\alpha$ (proof in Appendix).

\begin{thm}\label{thm_bdd} {\bf [MoF decreases with $\alpha$]}
There exists a $1\le B<\infty$ such that $\mof\le B^{1/\alpha}-1$ for every $\alpha>0$. \eop
\end{thm}
 From \eqref{eqn_choice_station}, 
 FoPS$(\alpha)$  with $\alpha>0$ would not allocate server to a bad-station (i.e., one with  bad travel condition) unless the latter has accumulated too little a utility by the decision epoch and hence are also opportunistic in nature; they attempt to maintain a balance between fairness and efficiency.  From  Theorem \ref{thm_bdd}, as $\alpha\to \infty$, the measure of fairness $\mof$ decreases to zero, implying almost identical utilities for all the stations; this is true even for the case with heterogeneous demands. However one now has to understand the PoF, and this is possible if probably the MoF of the efficient scheduler is analyzed. 
 
\noindent\textbf{Fairness of Efficient scheduler:} 
 An efficient scheduler maximizes system utility  and is known to be unfair  in the context of communication networks (e.g., \cite{cellular,mayurNH}).
 The systems appear to follow a different logic, when the schedulers are opportunistic; we have an interesting observation  about the fairness  of efficient scheduler itself.
  In \cite{mayurNH}, the authors argue that the PoF for opportunistic schedulers decreases with  increase in the number of agents and our results mirror the same.

The efficient scheduler usually attempts to schedule good stations, i.e., stations with good travel condition, and does not pay  attention to the stations that have accumulated low utility and have been starved.  However with $m$ sufficiently large, we show that there always exist  good non-empty stations. 
\begin{lemma}
Let ${\cal G}_k :=  \{ i: N_{i,k} > 0 , C_{i,k} \in \{g,0\} \}$. Then 
$P ( \mbox{support}(\bbeta_k^\alpha)  \subset  {\cal G}_k  )  \to 1$ as $m \to \infty$ for any $k$, when $\alpha =0$. 
\end{lemma}
{
\textbf{Proof:} 
At any decision epoch   consider, if possible, any $i \in {\cal G}$ and any $j \in {\cal G}^c$.    Then,

\vspace{-4mm}
{\small
\begin{eqnarray*}
O_i(s)-O_j(s)   &= &   w + \sum_{k=1}^{m}  \left ( { \lh}_k(j,s)-\lh_k(i,s)\right)   - { \gidc}_j(s), \\
&> & w  - { \gidc}_j(s)  \ge 0,
\end{eqnarray*}}as either ${ \gidc}_j(s)  = 0$  or $C_j = b$.
  Thus  any $j \notin {\cal G}$ is not scheduled, if the set ${\cal G}$ is non-empty. The probability, $P (|{\cal G}_k| = 0) \ge P \left ( \cap_{i\le m} \left \{ C_i = b,  {\cal N}_i (J) = 0 \right \} \right )  \to 0 $ as $m \to \infty$ for any $k.$  \eop
}

In view of the above Lemma, one can anticipate that the MoF/PoF of the system with large $m$ is close to that of an ideal system with only good travel conditions.  If further   $\lambda_i = \lambda$ for all $i$, then  $\mof$ of the  ideal system is zero -- from equation \eqref{eqn_fixed_point}, $\rsbu_i$ is the same for all $i$ by symmetry for any given $\alpha$.
Thus \textit{with large $m$, the MoF of efficient scheduler itself is small}. Further the ideal system has identical station utilities
 for all $\alpha$,  hence one can anticipate that the PoF is also zero.

\section{Numerical results}
In this section, we present few examples that illustrate
the variations in the measure of fairness $\mof$ and the price of fairness $\pof$ for different levels of fairness, number of stations, and with possibly heterogeneous demands. In the first example (Figure  \ref{fig:mof_vs_alpha}), travel towards two stations (from any starting point) is in bad   condition with probability 0.1, and for the remaining with probability 0.9; travel time follows Normal distribution (parameters in caption).  We observe that $\mof$ is  high for small $\alpha$, and then decreases with increase in $\alpha$, for any $m$ as discussed in Theorem \ref{thm_bdd}.  Further,  $\pof$  (right sub-figure) is smaller for higher $\alpha$, more importantly it is negligible for all the cases.

In another example in Figure \ref{fig:mof_and_util}, half the stations are good stations, i.e., have good travel conditions with high probability, and half are bad, further, the bad stations have lesser demand ($\lambda$). In the right sub-figure, we plot the average utility of good and bad stations as a function of $m$,  for different $\alpha$. When there are only two stations, with $\alpha=0$, the bad station has a negative utility; this goes in line with the fact that efficient scheduler prefers only good stations. With $\alpha>0$, both the stations have positive utility and disparity is lesser. As $m$ increases,  the utility of a bad station approaches that of a good station even for $\alpha=0$ (blue and red curves with only markers).\textit{ This reaffirms the fairness of efficient scheduler in a system with large $m$}. In the left sub-figure, we observe that $\mof$  is  smaller with higher $\alpha$ as in Theorem \ref{thm_bdd}; further, it also decreases with $m$. Interestingly, server  utility  (average of individual station utilities) is almost the same for $\alpha=0$ and $\alpha=1$; hence from \eqref{eqn_def_pof} the PoF is negligible.  However MoF $\mof$ (disparity in utilities) is higher,  for $m \ge 4$ and is close to zero only for $m > 10$.
Thus the PoF/MoF is  negligible even for the case with heterogeneous demands either with  large $\alpha $ or with large number of stations.

In all the figures, we have an interesting observation for $m=2,3$, the PoF   and the MoF  are both negligible. Thus fairness aspects are critical for systems with intermediate  number of stations.




{\bf Acknowledgement.} The Masters thesis of Mr Rishabh Kumar has inspired us to write this paper.

\section{Conclusion} 

We consider   polling systems with losses (due to customers lost),   gains (by serving) and opportunities (scheduler has some information about travel conditions). 
Assuming just the knowledge of binary indicator of the travel conditions, we   propose a family of fair schedulers that can achieve any desired level of fairness. Fairness is measured in terms of disparities in utilities of the individual stations. 

Using partial theoretical arguments and a conjecture, we illustrate the following: i) when the fairness parameter ($\alpha$) increases, the proposed scheduler achieves max-min fairness, i.e., the disparity in individual utilities reduces to zero; ii) as the number of stations increases, the efficient scheduler itself achieves the max-min fairness and the price of fairness reduces. The numerical examples illustrate that sometimes the schedulers are fair without paying much price, for number of stations as less as ten.

{\it Future work:} The work is incomplete without the proof of the conjecture. Further, the dynamic decisions in  queuing systems should also consider the future cost as in Markov Decision Process framework. A well posed sequential decision problem  which further caters to opportunistic   and fairness aspects would be an interesting topic for future work.

\hide{
We assume the knowledge of binary indicators of the travel conditions and propose a family of schedulers that cater to various levels of fairness. 


One of the major tasks in polling systems is server allocation.

An efficient scheduler maximizes the system utility,
and in wireless context the scheduler with $\alpha=0$ is known to be efficient (which maximizes the sum of the average utilities); this is true mainly because the utilities are i.i.d. in \eqref{eqn_beta_alpha_FP}. However FoPS scheduler with $\alpha=0$ need not be efficient in general, as one would actually require MDP framework to derive the efficient solution (see the discussions in footnote \ref{foot_one}); the MDP framework is required as the decision should consider not just the optimality of the instantaneous cost, but also that of the future cost. 

In this paper, we propose a set of schedulers which for the first time consider fairness aspect in  such a broad manner; basically the proposed schedulers can achieve various levels of fairness dictated by a single parameter. The schedulers may not be optimal in full generality, nonetheless consider dynamic decisions travel conditions and provide fairness.
}

\begin{appendix}

 { 
\textbf{Proof of Theorem \ref{Thm_conv}:} \textbf{Step 1.} Using law of large numbers and the fact that conditional utilities are bounded in either direction, i.e. (see \eqref{eqn_loss_bound}), 
$$-\lh^* \le \u_{i,k} \le w \mbox{ for all } i, k \mbox{ and sample paths } \omega,  $$
one can prove that the time averages $\sbu_{i,k}$ in \eqref{eqn_u_bar_average} are uniformly bounded between $[-  \lh^*, w]$, for all $k \ge K(\omega)$ and for almost all sample paths $\omega$, with $\lh^*$ as in \eqref{eqn_loss_bound}.

\textbf{Step 2.} From \eqref{eqn_u_bar_average} and \eqref{eqn_util}, again with $\Delta := \max\{\lh^*,w \}+w$, 
$$|\sbu_{i,k+m}-\sbu_{i,k}| 
\le \sum_{q=1}^{m} \frac{\Delta}{(k+q)^\gamma} \mbox{ for all } k\ge K(\omega),$$%
Thus $\{\sbu_{i,k}\}$ is a Cauchy sequence (a.s.) and hence converges to a  limit $\rsbu_i \in [-l^*,w]$  a.s., for every $i$. Observe that the limit $\bu = (\rsbu_1, \cdots, \rsbu_m)$  depends upon sample path $\omega$.

\textbf{Step 3.} From  \eqref{eqn_choice_station}, the transition probabilities of the  chain $\Y_k:= (S_{k-1}, \N_k, \C_k)$ freeze for all  $k> K(\omega)$ (redefine $K(\omega) $ if required to ensure this)  as $\bU_k \to \bu$; this is because the finitely many decisions in  \eqref{eqn_choice_station} remain the same for all $\bU_k \in \{\bu': |\bu' - \bu| \le \epsilon\}$ for some  appropriate $\epsilon > 0$. This step is true under {\bf  D}.0. 

\textbf{Step 4.} Thus $\{\Y_k\}$ is a homogeneous finite-state Markov chain for all $k\ge K(\omega)$.  We will have at least one irreducible closed class and  a stationary distribution (SD). Observe $\vta_k$ of \eqref{eqn_u_bar_act_average} equals time-average of an appropriate function of $\Y_k$; hence it converges to the  corresponding stationary expectation (depends upon $\Y_k$ at the time transition probabilities freeze)  by standard Ergodic theorems applicable for finite-state Markov chains. 
\eop

\textbf{Proof of Lemma \ref{lem_infinite_stay}:} Proof is by contradiction. If say under FoPS  the server stops visiting station $i$ after some $n$ epochs, i.e., say $S_k \ne i$ for all $k \ge n$. Then from \eqref{eqn_u_bar_average} and \eqref{eqn_actual_util}
there exist a finite $\tau$, such that $\bU_{i,t}<\delta$ and $N_{i,t}\ge 1$ for all $t\ge \tau$ with $\delta$ as in \eqref{eqn_dyn_dec}; so ${ \gidc}_i(\Z_t) = w$ for all $t\ge \tau$. 

Consider a $t\ge \tau$, for which $N_{S_{t-1},t}=0$ (the queue at current server location is empty), $C_{i,t}=g$; such a $t$ always exists because of the non-zero probability associated with all the relevant events. Consider $O_i(\Z_t)-O_j(\Z_t)$ where $j$ is any station that is visited i.o., obviously one such station always exists, (dropping notation $\Z_t$, see \eqref{eqn_ob} and say $S_{t-1}=s$):
\begin{equation*}
O_i(s)-O_j(s) =\frac{w}{\delta^\alpha} - \frac{{ \gidc}_j(s)}{\raU_j^\alpha}+\sum_{l} \frac{ { \lh}_l(j,s)-{ \lh}_l(i,s)}{\raU_l^\alpha}.\\
\end{equation*}When $j=s$ from \eqref{eqn_gains} we have ${ \gidc}_j(s) = 0$ and   then  (as $\raU_j \ge \delta$)
\begin{equation}\label{eqn_lem_1_proof}
  O_i(s)-O_j(s)\ge \frac{w-\sum_{l=1}^m  {\lh}_l(i,s)}{\delta^\alpha} > 0,
\end{equation}by assumption {\bf A}.1. For $j\ne s$, $ { \lh}_l(j,s)\ge { \lh}_l(i,s)$ since $C_{i,t}=g$, 
\begin{equation*}
  O_i(s)-O_j(s)\ge \frac{w}{\delta^\alpha} -  \frac{{ \gidc}_j(s)}{\raU_j^\alpha} \ge  0,
\end{equation*}
and then since the ties are broken randomly and equally likely, 
station $i$ would be chosen at some such $t$ from \eqref{eqn_choice_station}, thus a contradiction.  \eop

 \textbf{Proof of Theorem \ref{thm_bdd}:} We provide the proof by contradiction. For the given $\alpha$, say $(\rsbu_1,\dots,\rsbu_m)$ is the fixed point of \eqref{eqn_fixed_point}. Without loss of generality, assume that $\rsbu_1<,\dots,<\rsbu_m$ and if possible, say $\frac{\rsbu_m}{\rsbu_1}>B^{1/\alpha}$ where $B  > 1$ is such that, 
 $$B
> \frac{w}{
\min_{i,j,k} E[ {\cal N}_k (T_{i,j}) | C_{j} = g ]
 }.$$%
 %
%
We drop $\Z$ from the notation for simplicity and retain only $S$ and say, $S=s$ with $s\ne m$.  Consider $O_s(s)-O_m(s)$,

\vspace{-4mm}
{\small 
\begin{eqnarray*}
O_s(s)-O_m(s) \   \hspace{-25mm}\\
&&=\frac{{ \gidc}_s(s)-{ \lh}_s(s,s)+ { \lh}_s(m,s)}{\raU_s^\alpha} - \frac{{ \gidc}_m(s)-{ \lh}_m(m,s)+ { \lh}_m(s,s)}{\raU_m^\alpha}\\
&&\hspace{1cm}+\sum_{j\ne s, m} \frac{ { \lh}_j(m,s)-{ \lh}_j(s,s)}{\raU_j^\alpha},\\
&&=\frac{1}{\raU_m^\alpha}\left(\frac{\raU_m^\alpha}{\raU_s^\alpha}\gidc_s(s)- { \gidc}_m(s)+\sum_{j=1}^m \frac{\raU_m^\alpha}{\raU_j^\alpha}( { \lh}_j(m,s)-{ \lh}_j(s,s))\right),\\
%
&&\ge \frac{1}{\raU_m^\alpha}\left( \frac{\raU_m^\alpha}{\raU_1^\alpha}( { \lh}_1(m,s)-{ \lh}_1(s,s))
- { \gidc}_m(s)\right)>0,
\end{eqnarray*}} by the  choice of $B$. Hence,
$\beta_m(\Z) =0 $ if $S \ne m$.  This implies $P (S'= m | S \ne m) = 0$, which in turn implies that the state space is reducible under  limit transition probabilities,  with  class of states 
 $\{ S \ne m\}$ being closed.
 If $P(S= m) = 0$, then from \eqref{eqn_fixed_point},  $\rsbu_m < 0$ which contradicts $\raU_m > \raU_1 B \ge \delta B$. Thus $P(S= m) > 0$ (this is under SD), which implies $\{S = m\}$  is also closed class under  limit transition probabilities, and further that $P(S= m) = 1$; this also implies $P(S \ne m) = 0$. Hence under the SD, $N_i = b_i$ and $\raU_i = \delta$ for all $i < m$. Now the limiting chain only schedules $m$-th station, and thus  by {\bf A}.1 (and  standard  arguments for queues with load factor  less than one)   $P(N_m = 0) > 0$.
 But as argued in \eqref{eqn_lem_1_proof}, with  $N_1 >0$, $N_m = 0$ and $\raU_1 = \delta$ we have
 $O_1 \ge O_m$, which in turn implies $P(S = 1) >0$. This is again a contradiction and hence the assumption that $\mof > B^{1/\alpha}-1 $ is not true. \eop

}
\end{appendix}

\end{document}